\documentclass[11pt]{amsart}
\usepackage{amssymb,amsmath, graphicx, verbatim,epsfig,color}
\usepackage[latin1]{inputenc}
 \usepackage[all,cmtip]{xy}

\newcommand{\C}{\mbox{${\mathbb C}$}}

\newcommand{\p}{\mbox{${\mathbb P}$}}
\newcommand{\A}{\mbox{${\mathbb A}$}}

\newtheorem*{lm}{Lemma}

\newtheorem*{tm}{Theorem}

\theoremstyle{definition}\newtheorem*{rem}{Remark}

\newcounter{remarkcount}
\newcommand{\barr}{\overline}

\newcommand{\smap}{\barr{M}}

\newcommand{\st}{_{g,s}(X,\beta)}

\title{Gromov-Witten invariants and rational curves on Grassmannians}
\author{Alberto L\'opez Mart\'in}

\address{Institut f\"ur Mathematik, Universit\"at Z\"urich,
Z\"urich, CH-8057}
\email{alopez@math.uzh.ch}

\begin{document}
\bibliographystyle{plain}

\begin{abstract}
We study the enumerative significance of the $s$-pointed genus zero Gromov-Witten invariant on a homogeneous space $X$. For that, we give an interpretation in terms of rational curves on $X$. 
 \end{abstract}
\maketitle

\section{Introduction}

Since their appearance in the algebraic context, Gromov-Witten
invariants have proven to be an indispensable tool for enumerative
geometry. The problem of determining the number $N_d$ of rational
curves of degree $d$ passing through $3d-1$ points in general position
in the complex projective plane $\p^2$ was solved, by means of
Gromov-Witten theory, by Kontsevich (see\ \cite{KM94}).

Gromov-Witten invariants arose as enumerative invariants of stable
maps, which had been previously introduced independently by Ruan and
Tian \cite{RT95} in the symplectic case, and by Kontsevich and Manin
\cite{KM94} in the algebraic case. Let $X$ be a smooth projective
variety over $\C$, and let $\beta$ be a curve class on $X$. The set of
isomorphism classes of pointed maps $(C,p_1,\dots,p_s, f)$, where
$C$ is a projective nonsingular curve and $f$ is a morphism from $C$
to $X$ with $f_*([C])=\beta$, is denoted as $M\st$. Its compactification,
the moduli space $\smap\st$, parameterizes stable maps. The
stability condition is equivalent to finiteness of automorphisms of
the map. 

The purpose of our work is to study the enumerative significance of
genus zero Gromov-Witten invariants with a Grassmannian target.
While Gromov-Witten theory in the algebraic context in general requires
the sophisticated machinery of virtual fundamental classes
\cite{LT98, BF97} and the invariants have no clear
enumerative significance in general,
in the case
of genus zero with target a homogeneous variety $X=G/P$ the moduli stack is smooth
of the expected dimension and it makes sense to ask for the
enumerative significance of the Gromov-Witten invariants.
Given a tuple of classes of subvarieties of $X$ of suitable dimensions,
the Gromov-Witten invariant counts
isomorphism classes of $(C\cong\p^1,p_1,\ldots,p_s,f)$ such that
$f_*[\p^1]$ is equal to a given curve class and $f$ maps $p_i$ into
a general translate of the $i$th subvariety for all $i$,
according to Fulton and Pandharipande \cite{FP97}. (The $p_i$ are not fixed here; as explained in op. cit., page 93, there are alternative invariants in which the $p_i$ are fixed points in $\p^1$.) The purpose of this article is to rephrase this enumerative interpretation in
terms of rational curves on $X$ satisfying incidence conditions.
For a large class of enumerative conditions, we are able to exclude that due to reparameterizations of the source curve a map
$f$ contributes multiply to the Gromov-Witten invariant.
Then the Gromov-Witten invariant simply counts rational curves in
a given curve class incident to general translates of the given subvarieties.
The main Theorem asserts this to be the case when $X$ is a Grassmannian
variety and the subvarieties are Schubert varieties, up to a correction factor
of the degree of the curve class for each Schubert variety of codimension one.

\section{Preliminaries}
\label{prelim}

Let $G$ be a complex simple Lie group of classical type and $P$ a maximal parabolic
subgroup. The homogeneous space $X=G/P$ is a Grassmannian variety, a usual
Grassmannian of subspaces of some finite-dimensional complex vector space $V$
when $G$ is of type $A$, or of subspaces
isotropic for a given nondegenerate symmetric or
skew-symmetric bilinear form on $V$ in the other classical Lie types.
Throughout we assume that $\dim X\geq 2$.

The \emph{quantum Schubert calculus} is a set of
combinatorial rules that determine the genus zero three-point
Gromov-Witten invariants of $X$.
Quantum analogues of the classical Pieri and Giambelli formulas
are given for the usual Grassmannians by Bertram \cite{Ber97} and
for isotropic Grassmannians by
Buch, Kresch, and Tamvakis
\cite{KT03,KT04,BKT09,BKTpre,BKTinprep}.
In type $A$, an explicit combinatorial rule for the invariants, i.e.,
a quantum Littlewood-Richardson rule,
is available, due to Coskun \cite{Cos09}.

For the sake of completeness and to fix notation, let us recall some
definitions.
The Schubert varieties on the usual Grassmannian of $m$-planes in
$V\cong \C^n$ are indexed by integer partitions of length at most
$m$ and biggest part at most $n-m$.
There are analogous descriptions in the other Lie types,
based on \emph{$k$-strict partitions}; the reader is referred to
\cite{BKT09} for a detailed
description and facts about Schubert varieties in the
isotropic Grassmannians.
Here $k$ is $n-m$ when the bilinear form on $V$ is
skew symmetric and $X=\mbox{IG}(m,2n)$ and in the case of a symmetric bilinear form
$k$ is $n-m$, respectively $n+1-m$, when $X$ is
$\mbox{OG}(m,2n+1)$, respectively $\mbox{OG}(m,2n+2)$;
a partition is $k$-strict when there are no repeated parts $>k$.
We will generally denote by $X_\lambda$
the Schubert variety in $X$ corresponding to a partition $\lambda$.
Its codimension is $|\lambda|$, where for
$\lambda=(\lambda_1\geq\lambda_2\geq\dots\geq
\lambda_{\ell}>0)$ the weight is
$|\lambda|:=\lambda_1+\cdots+\lambda_\ell$.
The kernel and span of a rational map are
important notions, introduced and discussed by Buch in \cite{Bu03}. 
A rational map of degree $d$ to $X$ is a morphism
$f:\p^1\rightarrow X$ such that
$f_*[\p^1]$ has degree $d$.
Here degree is understood with respect to the projective embedding of
$X$ corresponding to a fundamental representation of $G$ with point
stabilizer $P$.
In dimension $1$, the unique Schubert class has degree $1$
and the corresponding Schubert variety is a line under this embedding.
Since $\dim X\geq 2, X$ contains a projective plane or a nonsingular quadric threefold, and therefore there exist rational curves
of every degree $d\geq 1$ on $X$. 

The moduli space of stable genus zero degree $d$ maps
$\overline{M}_{0,s}(X,d)$ is smooth (as a stack) of dimension
$\dim X+s-3+d\deg c_1(X)$.
There are evaluation maps $ev_i:\overline{M}_{0,s}(X,d)\rightarrow X$.
If $\alpha_1$, $\ldots$, $\alpha_s\in A^*(X)$ are classes in the
Chow (or cohomology) group of $X$ whose codimensions sum to
$\dim \overline{M}_{0,s}(X,d)$ then there is
the (non-gravitational, genus zero) \emph{Gromov-Witten invariant}
$$I_{d}(\alpha_1\cdots\alpha_s):=
\int_{\overline{M}_{0,s}(X,d)} ev_1^*\alpha_1\cup\cdots\cup ev_s^*\alpha_s.$$
If $\Gamma_1$, $\ldots$, $\Gamma_s$ are subvarieties of
$X$ with $\alpha_i$ (Poincar\'e dual to) the
fundamental class of $\Gamma_i$ for each $i$, then for
general $(g_1,\ldots,g_s)\in G^s$ the Gromov-Witten invariant
is equal to the number of degree $d$ maps
$\p^1\to X$ sending the $i$th marked point into $g_i\Gamma_i$ for each $i$.
For proofs of these facts, see \cite{FP97}.

\section{Main result}

We adopt the notation of Section \ref{prelim} and prove the following result.

\begin{tm}
If $X=G/P$ where $G$ is a complex simple Lie group of classical type and $P$ is a maximal
parabolic subgroup, $d$ and $s$ are positive integers,
and $\Gamma_1$,$\ldots$,$\Gamma_s$ are Schubert varieties whose
codimensions sum to $\dim \overline{M}_{0,s}(X,d)$,
then the Gromov-Witten invariant
$I_{d}([\Gamma_1]\cdots[\Gamma_s])$ is divisible by
$d^r$, where $r=\#\{i:\mathrm{codim}\,\Gamma_i=1\}$, and the quotient is equal to the
number of degree $d$ rational curves on $X$ incident to general translates
of the $\Gamma_i$.
\end{tm}

We note that if $\Gamma_i=X$ for any $i$ then the Gromov-Witten invariant is
zero (fundamental class axiom) and by Lemma 14 in \cite{FP97} the set of degree $d$ rational curves on $X$ incident to general translates
of the $\Gamma_i$ is empty. If some $r\geq 1$ of the $\Gamma_i$ have codimension 1,
then the Gromov-Witten invariant is equal
to $d^r$ times the $(s-r)$-point Gromov-Witten invariant with the
divisor classes omitted (divisor axiom).
This implies the divisibility assertion.
We first prove the enumerative claim assuming $\mathrm{codim}\,\Gamma_i\ge 2$ for
all $i$, then we
treat the case when some of the $\Gamma_i$ are divisors.

\begin{lm}
Let $\Gamma$ be a Schubert variety in $X$ of codimension at least $2$
with Schubert cell $\Gamma^0$.
Fix a point $p_1\in \p^1$.
Then for each $d\ge 1$ there exists a degree $d$ map
$f:\p^1\to X$ such that
\begin{enumerate}
\item $f$ is an unramified morphism;
\item $f(p_1)\in \Gamma^0$;
\item $f$ maps a nonzero tangent vector at $p_1$ to a tangent vector at
$f(p_1)$ not contained in the tangent space to $\Gamma^0$;
\item $f^{-1}(\Gamma)=\{p_1\}$.
\end{enumerate}
\end{lm}

Recall that $\overline{M}_{0,s}(X,d)$ is irreducible
\cite{Tho98, KP01}.
For a point $x\in X$ and for each $i$ we observe that by the group action there is a
birational isomorphism between $ev_i^{-1}(x)\times \A^{\dim X}$ and
$\overline{M}_{0,s}(X,d)$, and hence
$ev_i^{-1}(x)$ is irreducible as well.
In the situation of the Lemma we therefore obtain that
$ev_i^{-1}(\Gamma)$ is irreducible, by \cite[Theorem 4.17]{DM69}.

\begin{proof}
It suffices to verify (ii)-(iv) for a point of
$\overline{M}_{0,1}(X,d)$,
since the combination of these is
an open condition in $ev_1^{-1}(\Gamma)$ and (i) is
satisfied on a dense open subset of
$ev_1^{-1}(\Gamma)\cap M_{0,1}(X,d)$.
When $d=1$ it is clear that (ii)-(iv) may be satisfied.
When $d=2$ and $X$ is an orthogonal Grassmannian, this follows
from the description in
\cite[Lemma 2.1, proof of Thm.\ 2.3 or
Lemma 3.1, proof of Thm.\ 3.3]{BKT09}.

Otherwise, we consider two cases.
There is a critical degree $d_0$, the smallest for which two general points
on $X$ are joined by a rational curve of that degree:
$d_0=\min(m,n-m)$ when $X=G(m,n)$; otherwise $d_0$ is $m$ when $X=IG(m,2n)$
and $m$ rounded up to the next even integer for the orthogonal Grassmannians
(divided by two in the case of the maximal orthogonal Grassmannians).
The first case we consider is when $d\le d_0$.
Then we have the following set-up from
\cite[\S 2.2, 3.2, 4.2]{BKTJAMS},
\cite[\S 1.4, 2.4, 3.4]{BKT09}; specifically:
\begin{itemize}
\item a variety $Y_d$ parameterizing pairs $(A,B)$ with $\dim A=m-d$, $\dim B=m+d$, $A\subset B$, and $B\subset A^{\perp}$ when $X$ is an isotropic Grassmannian;
\item an incidence correspondence
$$\xymatrix{
T_d\ar[r] \ar[d]_\pi & X \\ Y_d
}$$
(where $T_d$ consists of triples $(A,\Sigma,B)$ with $A\subset\Sigma\subset B$ and $\Sigma$ a point of $X$);
\item ``modified'' Schubert varieties $Y_\lambda\subset Y_d$ each defined as the image by $\pi$ of the preimage $T_\lambda\subset T_d$ of $X_\lambda$;
and \item a result identifying the three-point genus zero Gromov-Witten invariants
in degree $d$ with (in some cases up to a certain power of $2$)
intersection points of modified Schubert varieties in $Y_d$.
(In types $B$ and $D$ when $d$ is odd there is a codimension $3$ subvariety
of $Y_d$ denoted $Z_d^\circ$ in \cite{BKT09} to which we need to
restrict our attention.) \end{itemize}
The fiber of $\pi$ above a general point of $Y_\lambda$ is a Grassmannian of
particular type (e.g., the Lagrangian Grassmannian of a symplectic $2d$-space
in type $C$).
Letting $\lambda^+$ be the smallest partition containing $\lambda$ so that
according to
\cite[proof of Cor.\ 2, 4 or 6]{BKTJAMS} or
\cite[Lemma 1.3, 2.1 or 3.1]{BKT09}
the map $T_{\lambda^+}\to Y_{\lambda^+}$ is generically finite,
the inequality $\dim Y_\lambda\geq \dim Y_{\lambda^+}=\dim T_{\lambda^+}$
is enough to guarantee that
the intersection of $T_\lambda$ with the fiber of $\pi$ above a general point
of $Y_\lambda$ has codimension at least $2$ in the fiber and is generically smooth.
Choose a general point of the intersection to be $f(p_1)$;
then a general rational map of degree $d$
in the fiber of $\pi$ meets requirements (ii)-(iv), since we have an explicit description of a general rational curve on the fiber of $\pi$ by \cite[Prop. 1.1]{BKT09}.
For $d>d_0$, we simply have to take a degree $d_0$ curve as just described
and attach $d-d_0$ copies of a degree $1$ tail.
If the point of attachment is general with respect to $f(p_1)$, then
a general line will for dimension reasons be disjoint from the
Schubert variety, and (ii)-(iv) remain valid.
\end{proof}

\begin{proof}[Proof of the Theorem]
In \cite[Lemma 14]{FP97} it is proven that, in the conditions of the Theorem, the intersection
$ev_1^{-1}(g_1\Gamma_1)\cap\cdots\cap ev_s^{-1}(g_s\Gamma_s)$
is a finite set of reduced points, each corresponding to an
irreducible source curve $\p^1$ with automorphism-free map to the
target variety $X$,
for general $(g_1,\ldots,g_s)\in G^s$.
The number of these is the Gromov-Witten invariant.

We prove the Theorem first in the case when each of the $\Gamma_i$ has
codimension at least $2$.
We claim, for general
$(g_1,\ldots,g_s)\in G^s$, that each of these finitely many intersection points
satisfies (i)-(iv) of the Lemma, with $\Gamma=g_1\Gamma_1$.
Given this, we may repeat the argument with the other $\Gamma_i$ in place of
$\Gamma_1$ and find for general
$(g_1,\ldots,g_s)\in G^s$ that for each of the finitely many
intersection points $(C\cong \p^1, p_1,\ldots, p_s,f)$
and each $i$, the point $p_i$ is the unique one on $C$ having
image contained in $g_i\Gamma_i$.

To prove the claim, by homogeneity we may fix $g_1=e$ and verify the
conditions for general $(g_2,\ldots,g_s)\in G^{s-1}$.
Let $ev_1^{-1}(\Gamma_1^0)^*$ denote the open subset of
$ev_1^{-1}(\Gamma_1^0)\cap M_{0,s}(X,d)$ satisfying
(i)-(iv) of the Lemma.
Now we apply the Kleiman-Bertini theorem to the diagram
$$\xymatrix{
&ev_1^{-1}(\Gamma_1^0)^* \ar[d]\\
\Gamma^0_2\times\cdots\times \Gamma^0_s \ar[r] & X^{s-1}
}
$$
with action of $G^{s-1}$.
Then, for general $(g_2,\ldots,g_s)$ the intersection is a finite set
of reduced points, and there is no contribution to the intersection
from points of $ev_1^{-1}(\Gamma_1^0)$ not in
$ev_1^{-1}(\Gamma_1^0)^*$.
The claim is verified.

For the case when some of the $\Gamma_i$ are divisors, we repeat the
above argument but taking
$ev_1^{-1}(\Gamma_1^0)^*$ to be defined by only conditions
(i)-(iii) of the Lemma.
Then for each of the surviving intersection points, the map
$f:\p^1\to X$ has properties (i)-(iii) at every preimage point of $\Gamma$.
It follows that there are $d$ distinct choices for a marked point on $\p^1$
mapping to $\Gamma_1$.
Hence there is a $d^r$-to-one correspondence (with $r$ as in the statement
of the Theorem) between intersection points in
$\overline{M}_{0,s}(X,d)$ and degree $d$ rational curves on $X$
satisfying the incidence conditions.
\end{proof}

\begin{rem}
The three-point genus zero Gromov-Witten invariants are those
that arise as structure constants in the small quantum cohomolgy ring;
they are those given by the quantum Schubert calculus.
For these, according to the Theorem, the Gromov-Witten invariant precisely counts
rational curves on $X$ except in one case:
when $d=2$ and one of the $\Gamma_i$ has codimension 1, then the
Gromov-Witten invariant is twice the number of conics satisfying the
incidence conditions.
This is only possible when $X$ is an orthogonal Grassmannian, reflecting
the presence of $q^2$ terms in the quantum Pieri formula for multiplication
with the codimension $1$ Schubert class (where $q$ is the formal
parameter of the quantum cohomology ring).
\end{rem}

\subsubsection*{Acknowledgements.} The author would like to thank A. Kresch, his thesis advisor, for his constant guidance and support. The author was supported in part by the Forschungskredit No. 57104401 from the Universit\"at Z\"urich.

\bibliography{gwgrass}

\end{document}